\documentclass[12pt]{article}

\title{
Sharp non-asymptotic oracle inequalities for nonparametric heteroscedastic regression models.
}

\author{L. Galtchouk
\thanks{
Department of Mathematics,
 Strasbourg University
 7, rue Rene Descartes,
 67084, Strasbourg, France, 
 e-mail: galtchou@math.u-strasbg.fr }
 \and 
S. Pergamenshchikov\thanks{
 Laboratoire de Math\'ematiques Raphael Salem,
 Avenue de l'Universit\'e, BP. 12,            
  Universit\'e de Rouen,                  
   F76801, Saint Etienne du Rouvray, Cedex France,
 e-mail: Serge.Pergamenchtchikov@univ-rouen.fr}
}

\date{}

\usepackage{srcltx}
\usepackage{amssymb}
\usepackage{amsfonts}
\usepackage{amsmath}
\usepackage{amsthm}
\usepackage{enumerate}
\usepackage{multicol}

\newtheorem{theorem}{Theorem}[section]

\newtheorem{lemma}[theorem]{Lemma}

\newtheorem{remark}{Remark}[section]

\newcommand\cA{{\cal A}}

\newcommand\cL{{\cal L}}
\newcommand\cB{{\cal B}}

\newcommand\cD{{\cal D}}

\newcommand\ve{\varepsilon}
\newcommand\ov{\overline}

\def\bbr{{\mathbb R}}

\def\text#1{\hbox{#1}}
\def\proof{{\noindent \bf Proof. }}
\def\endproof{\mbox{\ $\qed$}}

\def\E{{\bf E}}

\def\L{{\bf L}}
\newcommand{\wh}{\widehat}
\newcommand{\wt}{\widetilde}

\newcommand\Er{\mbox{Err}}

\def\Chi{{\bf 1}}
\def\d{\mathrm{d}}
\def\build #1_#2{\mathrel{\mathop{\kern 0pt #1}\limits_{#2}}}
\newcommand{\zs}[1]{{\mathchoice{#1}{#1}{\lower.25ex\hbox{$\scriptstyle#1$}}
{\lower0.25ex\hbox{$\scriptscriptstyle#1$}}}}
\numberwithin{equation}{section}

\begin{document}

\maketitle

\begin{abstract}
An adaptive nonparametric estimation procedure is constructed for 
 heteroscedastic regression when the noise variance depends on the unknown regression.
A non-asymptotic upper bound for a quadratic risk (oracle inequality) is obtained.
\end{abstract}
\vspace*{5mm}

\vspace*{5mm}

\bibliographystyle{plain}
\renewcommand{\columnseprule}{.1pt}
\newpage

\section{Introduction}\label{I}

Suppose we are given observations
$(y_j)_\zs{1\le j\le n}$ which obey the heteroscedastic
regression equation
\begin{equation}\label{I.1}
y_j\,=\,S(x_j)+\sigma_\zs{j}\xi_j\,,
\end{equation}
where design points $x_j=j/n$, $S(\cdot)$
 is an unknown function to be estimated,
$(\xi_j)_\zs{1\le j\le n}$ is a sequence of i.i.d. random variables,
$\E\xi_1=0\,,\quad \E\,\xi^2_1=1\,,\ (\sigma_\zs{j})_\zs{1\le j\le n}$
are unknown volatility coefficients, which may depend on design points and
on unknown regression function $S$. 
 
 The models of type \eqref{I.1} with $\sigma_\zs{j}=\sigma_\zs{j}(x_j)$ were introduced in 
\cite{AkVk} 
as a generalisation
of the nonparametric ANCOVA model of 
\cite{YoBo}.
 It should be noted that
heteroscedastic regressions with this type of volatility coefficients
 have been encountered in econometric studies, namely, in consumer budget
studies utilizing observations on individuals with diverse incomes and in
analyses of the investment behavior of firms of different sizes 
(see \cite{GoQu}).
For example,
for  consumer budget problems one  uses there (see p. 83)
 some parametric version of 
model \eqref{I.1}  with the volatility coefficient   defined as
\begin{equation}\label{I.1'}
\sigma^2_\zs{j}=c_\zs{0}+c_\zs{1}x_\zs{j}+c_\zs{2}S^2(x_\zs{j})\,,
\end{equation}
where $c_\zs{0}$, $c_\zs{1}$ and $c_\zs{2}$ are some nonnegative  unknown constants.

 Moreover, this
regression model appears in the drift estimation problem for stochastic differential
equations when one passes from continuous time to discrete time model by making use of 
sequential kernel estimators having asymptotically minimal variances 
(see \cite{GaPe1},\cite{GaPe3}-\cite{GaPe5}).

The volatility coefficient estimation in heteroscedastic regression was considered in a few papers (see, for example, \cite{CaWa} and the references therein).
By making use of the squared
first-order differences of the observations the initial problem in that paper was reduced to the
regression function estimation in the model of type \eqref{I.1}. 

In this paper we develop the approach proposed in \cite{GaPe2}.
The first goal of the research is to construct an adaptive 
 procedure for estimating  the function $S$ which does not use any smoothness information of $S$ and
which is based on observations $(y_\zs{j})_\zs{1\le j\le n}$  and further
to obtain a sharp non-asymptotic upper bound (oracle inequality) for a quadratic risk in the case
 when the smoothness of $S$ is unknown. The second goal is to prove that the
 constructed procedure is efficient also in the asymptotic setup.


Problems of constructing a nonparametric estimator and proving a 
non-asymptotic  upper bound for a risk in homoscedastic model, that is when $\sigma_\zs{j}\equiv\sigma$,
 were studied in  few papers.
  A non-asymptotic upper bound  for a quadratic risk over thresholding estimators
 is given in 
 \cite{KaMa}.
  In papers 
\cite{BaBiMa}, \cite{Ma}
an adaptive model selection procedure has been constructed. It is
based on least squares estimators
 and a non-asymptotic upper bound has been obtained
for a quadratic risk which is best in the principal term for the given class of
estimators when the noise vector $(\xi_1\,\ldots,\xi_n)$ is gaussian.
 This type of upper bounds is called the {\sl oracle inequality}.
In 
\cite{FoPe} 
the oracle inequality has been obtained
for a model selection procedure based on any estimators 
in the case when
the noise vector $(\xi_1,\ldots,\xi_n)$ has a spherically symmetric distribution.
Moreover, 
some sharp oracle inequalities have been obtained also for homoscedastic
regression with gaussian noises, see, for example, 
\cite{Kn}.
 Here the adjective "sharp" means
that the coefficient of the principal term may be chosen as close to unity as desired.

In the paper for  heteroscedastic regression  an adaptive procedure
is constructed for which the sharp non-asymptotic oracle inequality is proved.
It should be noted that the methods used in former papers to obtain the sharp oracle inequality
in regression models are limited by the homoscedastic case since they are based
 on the fact that an orthogonal transformation of a noise gaussian vector
 $(\xi_1,\ldots,\xi_n)$ gives
 a gaussian vector. In heteroscedastic regression models under consideration these methods are not
 valid since the noise vector is not gaussian. To obtain  sharp non-asymptotic oracle inequalities
 in the heteroscedastic case 
the authors develop a new  mathematical
 tools based on "penalty" methods and  Pinsker's type weights.

Moreover, 
in \cite{GaPe6}
we show that the given adaptive estimator is 
efficient in the asymptotic sense, that is,
the sharp asymptotic lower bound is proved for a quadratic risk and it is attained
over this estimator.
The sharp non-asymptotic oracle inequality plays a cornerstone role in proving the asymptotic
efficiency. To obtain the optimal upper bound for the risk one should use a weighted least squares
estimator with weights depending on the smoothness of the unknown regression function. The smoothness
being unknown, one can't use directly this weighted least squares estimator giving the minimal upper 
bound. The given sharp non-asymptotic oracle inequality allow us to replace this unknown
weighted least squares estimator with an adaptive estimator. The risk of the adaptive estimator is
less than the risk of optimal (unknown) weighted least squares estimator up to additive and  
multiplicative constants. Taking in account that the multiplicative constant tends to one and the 
order of the additive constant is less then the order of the convergence rate, we obtain that the risk
of the given adaptive procedure asymptotically coincides with the risk of the
optimal (unknown) weighted least squares estimator. Therefore, given the optimal lower bound, we obtain
the asymptotic efficiency of the adaptive procedure satisfying the  sharp non-asymptotic oracle 
inequality.



    The paper is organized as follows.
In Section 2 we construct an adaptive estimation procedure based on
weighted least squares 
estimators and
 we obtain a non-asymptotic upper bound for the quadratic risk.
 In Section 3 we propose an estimator for the integrated noise variance and give 
 the oracle inequality in the case of Sobolev space, $S\in W^{k}_\zs{r}$.
The proofs are given in Section 4. Section 5 contains a numerical comparison of the given procedure
with an adaptive procedure proposed in \cite{EfPi}.
The Appendix contains some technical results.

\section{Oracle inequality}\label{N}

In this paper we study the non-asymptotic estimation problem of the function $S$
in the model \eqref{I.1} by observations $(y_\zs{j})_\zs{1\le j\le n}$
 with odd sample number $n$.
 We assume that in \eqref{I.1}
 the sequence
$(\xi_j)_\zs{1\le j\le n}$
is i.i.d. with
\begin{equation}\label{Or.1}
\E\xi_1=0\,,\quad
\E\,\xi^2_1=1
\quad\mbox{and}\quad
\E\xi^4_1=\xi^*<\infty\,.
\end{equation}

Moreover, we assume that $(\sigma_\zs{l})_\zs{1\le l\le n}$ is a sequence of positive
random variables independent
of $(\xi_i)_\zs{1\le i\le n}$ and bounded away from $+\infty$, i.e.
there exists some nonrandom unknown constant $\sigma_*> 0$ such that
\begin{equation}\label{Or.2}
\max_\zs{1\le l\le n}\,\sigma^2_\zs{l}\,\le\,\sigma_*\,.
\end{equation}

For any estimate $\wh{S}_n$ of $S$ based on observations
$(y_j)_\zs{1\le j\le n}$, the estimation accuracy is measured by the mean integrated
squared error (MISE)
\begin{equation}\label{Or.3}
\E_\zs{S}\,\|\wh{S}_n-S\|_n^2\,,
\end{equation}
where
$$
\|\wh{S}_\zs{n}-S\|^2_n=(\wh{S}_\zs{n}-S,\wh{S}_\zs{n}-S)_n=
\frac{1}{n}\sum^n_\zs{l=1}(\wh{S}_\zs{n}(x_l)-S(x_l))^2\,.
$$

 We make use of the trigonometric
basis $(\phi_j)_\zs{j\ge 1}$ in $\cL_2[0,1]$ with
\begin{equation}\label{Or.4}
\phi_1=1\,,\quad
\phi_\zs{j}(x)=\sqrt{2}\,Tr_\zs{j}(2\pi [j/2]x)\,,\ j\ge 2\,,
\end{equation}
where the function $Tr_\zs{j}(x)=\cos(x)$ for even $j$ and
$Tr_\zs{j}(x)=\sin(x)$ for odd $j$; $[x]$ denotes the integer part of $x$.
Note that if $n$ is odd, then this basis is orthonormal for the empirical inner product generated by the sieve
$(x_\zs{j})_\zs{1\le j\le n}$, that is for any $1\le i,j\le n$,
\begin{equation}\label{Or.5}
(\phi_i\,,\,\phi_j)_\zs{n}=
\frac{1}{n}\sum^n_\zs{l=1}\phi_i(x_l)\phi_j(x_l)=\delta_\zs{ij}\,,
\end{equation}
where $\delta_\zs{ij}$ is Kronecker's symbol.

By making use  of this basis we define the discrete Fourier transformation in \eqref{I.1}
and obtain the Fourier coefficients
\begin{equation}\label{Or.6}
\wh{\theta}_\zs{j,n}=(Y,\phi_j)_n
\quad\mbox{and}\quad
\theta_\zs{j,n}=(S,\phi_j)_n\,.
\end{equation}
Here  $Y=(y_\zs{1},\ldots,y_\zs{n})'$ and
$S=(S(x_\zs{1}),\ldots,S(x_\zs{n}))'$. The prime denotes the transposition.

From \eqref{I.1} it follows directly that these Fourier coefficients satisfy the following
equation
\begin{equation}\label{Or.7}
\wh{\theta}_\zs{j,n}=\theta_\zs{j,n}+\frac{1}{\sqrt{n}}\xi_\zs{j,n}
\quad
\mbox{with}
\quad
\xi_\zs{j,n}=\frac{1}{\sqrt{n}}\sum^n_\zs{l=1}\sigma_\zs{l}\xi_\zs{l}\phi_j(x_l)\,.
\end{equation}

We estimate the function $S$ by the weighted least squares estimator
\begin{equation}\label{Or.8}
\wh{S}_\zs{\lambda}(x)=\sum^n_\zs{j=1}\lambda(j)\wh{\theta}_\zs{j,n}\phi_\zs{j}(x)\,,
\end{equation}
where $x\in[0,1]$, the weight vector $\lambda=(\lambda(1),\ldots,\lambda(n))'$
belongs to some finite set $\Lambda$ from $[0,1]^n$. We denote by $\nu$ the cardinal number of the set $\Lambda$.

Now we need to write a cost function to choose a weight $\lambda\in\Lambda$. Of course,
it is obvious, that the best way is to minimize the cost function which is equal to 
the empirical squared error
$$
\Er_\zs{n}(\lambda)=\|\wh{S}_\zs{\lambda}-S\|^2_\zs{n}\,,
$$
 which in our case is equal to
\begin{equation}\label{Or.9}
\Er_\zs{n}(\lambda)\,=\,
\sum^n_\zs{j=1}\,\lambda^2(j)\wh{\theta}^2_\zs{j,n}\,-
2\,\sum^n_\zs{j=1}\,\lambda(j)\wh{\theta}_\zs{j,n}\,\theta_\zs{j,n}\,+\,
\sum^n_\zs{j=1}\theta^2_\zs{j,n}\,.
\end{equation}
Since coefficients $\theta_\zs{j,n}$ are unknown, we need to replace the term
$\wh{\theta}_\zs{j,n}\,\theta_\zs{j,n}$
by some estimator which we choose as
$$
\wt{\theta}_\zs{j,n}=
\wh{\theta}^2_\zs{j,n}-\frac{1}{n}\wh{\varsigma}_\zs{n}\,,
$$
 where $\wh{\varsigma}_\zs{n}$  is some estimator of the integrated noise variance
\begin{equation}\label{Or.10}
\varsigma_\zs{n}=n^{-1}\,\sum^n_\zs{l=1}\,\sigma^2_\zs{l}\,.
\end{equation}
Such type of estimators is given in \eqref{Si.4}.

Moreover, for this substitution to the empirical squared error one needs to pay
a penalty. Finally, we define the cost function by the following way
\begin{equation}\label{Or.11}
J_\zs{n}(\lambda)\,=\,\sum^n_\zs{j=1}\,\lambda^2(j)\wh{\theta}^2_\zs{j,n}\,-
2\,\sum^n_\zs{j=1}\,\lambda(j)\,\wt{\theta}_\zs{j,n}\,
+\,\rho \wh{P}_\zs{n}(\lambda)\,,
\end{equation}
where $\rho$ is some positive coefficient which will be chosen later.
 The penalty term we define as
\begin{equation}\label{Or.12}
\wh{P}_\zs{n}(\lambda)=\frac{|\lambda|^2 \wh{\varsigma}_\zs{n}}{n}
\quad\mbox{with}\quad
|\lambda|^2=\sum^n_\zs{j=1} \lambda^2(j)\,.
\end{equation}
Note that in the case when the sequence $(\sigma_\zs{l})_\zs{1\le l\le n}$ is known,
i.e.  $\wh{\varsigma}_\zs{n}=\varsigma_\zs{n}$, we obtain
\begin{equation}\label{Or.13}
{P}_\zs{n}(\lambda)=\frac{|\lambda|^2 \varsigma_\zs{n}}{n}\,.
\end{equation}
We set
\begin{equation}\label{Or.14}
\wh{\lambda}=\mbox{argmin}_\zs{\lambda\in\Lambda}\,J_n(\lambda)
\end{equation}
and  define an estimator of $S$ as
\begin{equation}\label{Or.15}
\wh{S}_\zs{*}=\wh{S}_\zs{\wh{\lambda}}\,.
\end{equation}
We recall that the set $\Lambda$ is finite so $\wh{\lambda}$ exists. In the case 
when $\wh{\lambda}$ is not unique we take one of them.

To formulate the oracle inequality we introduce,
for $0<\rho< 1/3$,
 the following function
\begin{equation}\label{Or.16}
\Upsilon^*_\zs{n}(\rho)=\frac{16\nu}{\rho}
+4u_\zs{1,n}\left(1+\nu\frac{\sqrt{\xi^*}}{\sqrt{n}}\right)+
4\nu v_\zs{n}\frac{\sqrt{\xi^*}}{\sqrt{n}}\,.
\end{equation}
 Here and thereafter we make use of the following notations: for $i=1,2$,
\begin{equation}\label{Or.17}
v_\zs{n}=\max_\zs{\lambda\in\Lambda}\,\sum^n_\zs{j=1}\lambda(j)
\quad
\mbox{and}
\quad
u_\zs{i,n}=\max_\zs{\lambda\in\Lambda}\,\sup_\zs{1\le l\le n}
|\sum^n_\zs{j=1}\lambda^i(j)\,(\phi^2_\zs{j}(x_\zs{l})-1)|\,,
\end{equation}.

\medskip

\begin{theorem}\label{Th.Or.1}
Let $\Lambda$ be any finite set in $[0,1]^n$. For any $n\ge 3$ and\\  $0<\rho<1/3$,
the estimator $\wh{S}_\zs{*}$ satisfies the oracle inequality
\begin{align}\label{Or.18}
\E_\zs{S}\|\wh{S}_\zs{*}-S\|^2_\zs{n}\,
\le\,
\frac{1+3\rho-2\rho^2}{1-3\rho}
\,
\min_\zs{\lambda\in\Lambda}\,&
\E_\zs{S}\|\wh{S}_\zs{\lambda}-S\|^2_\zs{n}
+\frac{1}{n}\,\cB_\zs{n}(\rho)\,,
\end{align}
where
$\cB_\zs{n}(\rho)=\Psi_\zs{n}(\rho)+\kappa(\rho)v_\zs{n}\,\E_\zs{S}|\wh{\varsigma}_\zs{n}-\varsigma_\zs{n}|$
with
$$
\Psi_\zs{n}(\rho)=
\frac{
\rho(1-\rho)\Upsilon^*_\zs{n}(\rho)+2\nu
+2\rho^2(1-\rho)u_\zs{2,n}
}{\rho(1-3\rho)}\sigma_*
$$
and $\kappa(\rho)=4(1-\rho^2)/(1-3\rho)$.

If in the model \eqref{I.1}
the volatility coefficients $(\sigma_\zs{l})_\zs{1\le l\le n}$ are known,
 then\\
$\wh{\varsigma}_\zs{n}=\varsigma_\zs{n}$ and inequality \eqref{Or.18}
 has the following form
\begin{equation}\label{Or.19}
\E_\zs{S}\|\wh{S}_\zs{*}-S\|^2_\zs{n}\,
\le\,
\frac{1+3\rho-2\rho^2}{1-3\rho}\,
\min_\zs{\lambda\in\Lambda}
\E_\zs{S}\|\wh{S}_\zs{\lambda}-S\|^2_\zs{n}
+\frac{1}{n}\Psi_\zs{n}(\rho)\,.
\end{equation}

\end{theorem}

\medskip
\medskip

\begin{remark}\label{Re.3.1}
Note that the principal term in the right-hand side of \eqref{Or.18}--\eqref{Or.19} is best
in the class of estimators $(\wh{S}_\zs{\lambda}\,,\,\lambda\in\Lambda)$.
Inequalities of such type are called {\rm the sharp non-asymptotic
oracle inequalities}. The inequality is sharp in the sense that the coefficient
of the principal term may be chosen as close to $1$ as desired. Similar inequalities
for homoscedastic models \eqref{I.1} with $\sigma_\zs{l}=\sigma$ were given, for example, 
 in \cite{Kn}.
The methods used there cannot be extended to the heteroscedastic case since,
after the Fourier transformation, the random variables $(\xi_\zs{i,n})$ in model \eqref{Or.7}
are dependent contrary to the homoscedastic case, where these random variables are
independent (see, for example,\cite{Ro}). 
\end{remark}

\medskip

\begin{remark}\label{Re.3.2}
If one would like to obtain the asymptotically minimal MISE of the
estimator $\wh{S}_\zs{*}$, then the secondary term $\cB_\zs{n}(\rho)$
in \eqref{Or.18} should be slowly varing, i.e. 
for any $\gamma>0$,
\begin{equation}\label{Or.19-1}
\cB_\zs{n}(\rho)/n^{\gamma}\to 0\,,
\quad\mbox{as}\quad n\to\infty\,.
\end{equation}
 Indeed, since usually the optimal rate
is of order $n^{2k/(2k+1)}$ for some $k\ge 1$, then after multiplying the inequality
\eqref{Or.18} by this rate the principal term gives the optimal constant and the secondary
one is of the type 
$$
\cB_\zs{n}(\rho)/n^{1/(2k+1)}\,.
$$
Therefore the property \eqref{Or.19-1} provides the asymptotic vanishing,
as $n\to \infty$, the secondary term $\cB_\zs{n}(\rho)$ for $k\ge 1$. 
To obtain the property \eqref{Or.19-1}, it suffices that, for any $\gamma>0$,
$$
\rho n^{\gamma}\to +\infty\,,
$$
and
\begin{equation}\label{Or.19-2}
\frac{u_\zs{1,n}+ u_\zs{2,n}+v_\zs{n}\E_\zs{S}|\wh{\varsigma}_\zs{n}-\varsigma_\zs{n}|}{n^{\gamma}}\to 0\,,
\quad\mbox{as}\quad n\to\infty\,,
\end{equation}
 thanks to definitions of $\cB_\zs{n}(\rho)$ and $\Psi_\zs{n}(\rho)$.
To obtain the first convergence it suffices to take the parameter $\rho$ as $\rho \ge \varrho_\zs{n}$,
where $\varrho_\zs{n}$ is a slowly decreasing function, i.e.
\begin{equation}\label{Or.20}
\lim_\zs{n\to\infty}\,\varrho_\zs{n}\,=\,0
\quad
\mbox{and for any}
\quad \gamma>0\quad
\lim_\zs{n\to\infty}\, n^\gamma\,\varrho_\zs{n}\,=\,+\infty\,,
\end{equation}
for example, $\varrho=1/\ln n$. 
For the second convergence the choice of $u_\zs{1,n}, u_\zs{2,n}, v_\zs{n}$ and of the estimator 
$\wh{\varsigma}_\zs{n}$ is proposed below.
\end{remark}

\medskip

Consider now the order of the termes $v_\zs{n}, u_\zs{1,n}, u_\zs{2,n}$
and the function $\Psi_\zs{n}(\rho)$ in the case when the finite set $\Lambda$ is formed by
a special version of Pinsker's weights  (see, for example, \cite{Nu}). 
 To this end,
 we define the sieve
$$
\cA_\zs{\ve}=\{1,\ldots,k^*\}\times\{t_1,\ldots,t_m\}\,,
$$
where  $t_i=i\ve$ and $m=[1/\ve^2]$. We suppose that the parameters 
$k^*\ge 1$ and $0<\ve\le 1$ are functions of $n$, i.e. $k^*= k^*_\zs{n}$ and 
$\ve=\ve_\zs{n}$, such that,
\begin{equation}\label{Or.21}
\left\{
\begin{array}{ll}
&\lim_\zs{n\to\infty}\,k^*_\zs{n}=+\infty\,,
\quad
\lim_\zs{n\to\infty}\,k^*_\zs{n}/\ln n=0\,,\\[4mm]
&\lim_\zs{n\to\infty}\,\ve_\zs{n}\,=\,0
\quad\mbox{and}\quad
\lim_\zs{n\to\infty}\,n^{\gamma}\,\ve_\zs{n}\,=+\infty\,,
\end{array}
\right.
\end{equation}
for any $\gamma>0$. For example, one can take for  $n\ge 3$
\begin{equation}\label{Or.22}
\ve_\zs{n}=1/\ln n
\quad\mbox{and}\quad
k^*_\zs{n}=\ov{k}+\sqrt{\ln n}\,,
\end{equation}
where $\ov{k}$ is any nonnegative constant.

For any $\alpha=(\beta,t)\in\cA_\zs{\ve}$ we define the weight vector
$\lambda_\zs{\alpha}=(\lambda_\zs{\alpha}(1),\ldots,\lambda_\zs{\alpha}(n))'$ as
\begin{equation}\label{Or.23}
\lambda_\zs{\alpha}(j)=\Chi_\zs{\{1\le j\le j_\zs{0}\}}+
\left(1-(j/\omega_\alpha)^\beta\right)\,
\Chi_\zs{\{ j_\zs{0}<j\le \omega_\alpha\}}\,.
\end{equation}
Here $j_0=j_0(\alpha)=\left[\omega_\zs{\alpha}/\ln n\right]$ with
$$
\omega_\zs{\alpha}=\ov{\omega}+(A_\zs{\beta}\,t\,n)^{1/(2\beta+1)}\,,
$$
where $\ov{\omega}$ is any nonnegative constant and
$$
A_\zs{\beta}=(\beta+1)(2\beta+1)/(\pi^{2\beta}\beta)\,.
$$
Hence,
\begin{equation}\label{Or.24}
\Lambda\,=\,\{\lambda_\zs{\alpha}\,,\,\alpha\in\cA_\zs{\ve}\}
\end{equation}
and $\nu=\nu_\zs{n}=k^*_\zs{n} m_\zs{n}$. Note that in this case, in view of 
\eqref{Or.21}, for any $\gamma>0$
$$
\lim_\zs{n\to\infty}\,\nu_\zs{n}/n^{\gamma}=0\,.
$$
Moreover, by \eqref{Or.23}
\begin{align*}
\sum^n_\zs{j=1}\lambda_\zs{\alpha}(j)=
\Chi_\zs{\{j_\zs{0}\ge 1\}}\,j_\zs{0}
+\Chi_\zs{\{\omega_\zs{\alpha}\ge 1\}}
\sum^{[\omega_\zs{\alpha}]}_\zs{j=j_\zs{0}+1}
\left(1-(j/\omega_\zs{\alpha})^\beta\right)\le \omega_\zs{\alpha}\,.
\end{align*}
Therefore, taking into account that
$A_\zs{\beta}\le A_1<1$ for $\beta\ge 1$
we find that
$$
v_\zs{n}\le\ov{\omega}+(n/\ve )^{1/3}\,,
$$
i.e. 
\begin{equation}\label{Or.25}
\limsup_\zs{n\to\infty}\frac{v_\zs{n}}{\sqrt{n}}\,<\,\infty\,.
\end{equation}
Moreover, note that for any $1\le l\le n$, we get
\begin{align*}
\sum^n_\zs{j=1}\lambda_\zs{\alpha}(j)\,(\phi^2_\zs{j}(x_\zs{l})-1)&=
\Chi_\zs{\{j_\zs{0}\ge 1\}}
\sum^{j_\zs{0}}_\zs{j=1}\,
(\phi^2_\zs{j}(x_\zs{l})-1)\\
&+\Chi_\zs{\{\omega_\zs{\alpha}\ge 1\}}
\sum^{[\omega_\zs{\alpha}]}_\zs{j=j_\zs{0}+1}
\left(1-(j/\omega_\zs{\alpha})^\beta\right)
(\phi^2_\zs{j}(x_\zs{l})-1)\,.
\end{align*}
Thus
Lemma~\ref{Le.A.3}
implies that
$$
u_\zs{1,n}\le\,1+2^{\beta+1}\le
1+2^{ k^*+1}\,.
$$
Due to the condition for $k^*_\zs{n}$ in \eqref{Or.21}
 this function is slowly varying, i.e.
for any $\gamma>0$,
\begin{equation}\label{Or.26}
\lim_\zs{n\to\infty}
u_\zs{1,n}/n^\gamma=0\,.
\end{equation}
By the same way we obtain that
$$
u_\zs{2,n}\le
1+2^{ k^*+2}+2^{2k^*+1}
$$
and, therefore, for any $\gamma>0$
\begin{equation}\label{Or.27}
\lim_\zs{n\to\infty}
u_\zs{2,n}/n^\gamma=0\,.
\end{equation}
Then for any sequence  $(\varrho_\zs{n})_\zs{n\ge 1}$ satisfying properties \eqref{Or.20} and
for any $\gamma>0$,
$$
\lim_\zs{n\to\infty}\,\sup_\zs{\varrho_\zs{n}\le \rho\le 1/6}
\Psi_\zs{n}(\rho)/n^\gamma=0\,.
$$

\medskip
\medskip

\begin{remark}\label{Re.3.3}
As we shall see in the proof of Theorem~\ref{Th.Or.1}, the oracle inequality is true for any basis and 
any design $(x_\zs{k})_\zs{1\le k\le n}$ which possesse the orthonormality property \eqref{Or.5}.
Moreover, if the sequences \eqref{Or.17} and an estimator of the unknown integrated variance 
$\wh{\varsigma}_\zs{n}$ satisfy the property \eqref{Or.19-2}, then the secondary term in the inequality
\eqref{Or.18} possesses the property \eqref{Or.19-1}. In the next section we give an estimator
$\wh{\varsigma}_\zs{n}$ in the case of the trigonometric basis and the equidistant design.
 \end{remark}

\medskip

\section{Oracle inequality for $S\in W^{k}_\zs{r}$}

Assume that $S\,:\,\bbr\to\bbr$ is a $k$ times differentiable 1-periodic function
such that
\begin{equation}\label{Si.1}
\sum_\zs{j=0}^k\,\|S^{(j)}\|^2\le r\,,
\end{equation}
where
\begin{equation}\label{Si.2}
\|f\|^2=\int^1_\zs{0}f^2(t)\d t\,.
\end{equation}
We denote by $W^{k}_\zs{r}$ the set of all such functions.
Moreover, we suppose that $r>0$ and $k\ge 1$ are unknown parameters.

Note that,  the space $W^{k}_\zs{r}$ can be represented as
an ellipses in the Hilbert space, i.e.
\begin{equation}\label{Si.3}
W^{k}_\zs{r}=\{S\in\cL_\zs{2}[0,1]\,:\,S=\sum^\infty_\zs{j=1}\theta_\zs{j}\phi_\zs{j}
\quad\mbox{such that}\quad
\sum_\zs{j=1}^\infty\,a_\zs{j}\theta^2_\zs{j}\le r\}\,,
 \end{equation}
where the basis functions $(\phi_\zs{j})_\zs{j\ge 1}$ are defined in \eqref{Or.4};
$(\theta_\zs{j})_\zs{j\ge 1}$ are the Fourier coefficients, i.e.
\begin{equation}\label{Si.3-1}
\theta_\zs{j}=(S,\phi_\zs{j})=\int^1_\zs{0}S(t)\phi_\zs{j}(t)\d t\,.
\end{equation}
The coefficients $(a_\zs{j})_\zs{j\ge 1}$ are defined as
$$
a_\zs{j}=\sum^k_\zs{l=0}\|\phi^{(l)}_\zs{j}\|^2=
\sum^k_\zs{l=0}(2\pi [j/2])^{2l}\,.
$$
To estimate $\varsigma_\zs{n}$, we make use of the following estimator:
\begin{equation}\label{Si.4}
\wh{\varsigma}_\zs{n}=\sum^{n}_\zs{j=d_\zs{n}+1}
\wh{\theta}^2_\zs{j,n}\,,
\end{equation}
where the parameter
$d_\zs{n}\,,\ 1\le d_\zs{n}\le n-1\,,$ will be chosen later.

In Section~\ref{Se.P} we show the following result.
\begin{lemma}\label{Le.Si.1}
For any $n\ge 2$, $r>0$ and $S\in\,W^1_\zs{r}$
\begin{equation}\label{Si.5}
\E_\zs{S}\,|\wh{\varsigma}_\zs{n}\,-\,\varsigma_\zs{n}|
\le\frac{2\left(\sqrt{\xi^*}+\sqrt{2}\right)\sigma_*
+ \varsigma^*_\zs{n}(r)}{\sqrt{n}}\,,
\end{equation}
where 
$$
\varsigma^*_\zs{n}(r)=\frac{4r\sqrt{n}}{d^2_\zs{n}}
+4\sqrt{r\sigma_*}\frac{1}{d_\zs{n}}
+
\frac{(2+d_\zs{n})\sigma_*}{\sqrt{n}}\,.
$$
\end{lemma}
\noindent If we choose the parameter $d_\zs{n}$
in \eqref{Si.4} such that
\begin{equation}\label{Si.5-2}
\lim_\zs{n\to\infty}d_\zs{n}/\sqrt{n}=0
\quad\mbox{and}\quad
\lim_\zs{n\to\infty}d^2_\zs{n}/\sqrt{n}=\infty\,,
\end{equation}
 we obtain that 
$$
\lim_\zs{n\to\infty}\varsigma^*_\zs{n}(r)=0\,.
$$ 
Theorem~\ref{Th.Or.1} and  inequality \eqref{Si.5} imply immediately the
following result.
\begin{theorem}\label{Th.Si.1}
Let $\Lambda$ be any finite set in $[0,1]^n$. Assume that in the model \eqref{I.1} the function
$S$ belongs to  $W_\zs{r}^{1}$.
Then, for any $n\ge 3$ and $0<\rho<1/3$,
 the procedure $\wh{S}_\zs{*}$ from \eqref{Or.15} with $\wh{\varsigma}_\zs{n}$
defined by \eqref{Si.4} and \eqref{Si.5-2}
 satisfies the following oracle inequality
\begin{equation}\label{Si.6}
\E_S\|\wh{S}_\zs{*}-S\|^2_n
\le\,\frac{1+3\rho-2\rho^2}{1-3\rho}\,
\min_\zs{\lambda\in\Lambda}\E_S\|\wh{S}_\zs{\lambda}-S\|^2_n
+\frac{1}{n}\cD_\zs{n}(\rho,r)\,,
\end{equation}
where
$$
\cD_\zs{n}(\rho,r)=\Psi_\zs{n}(\rho)+\kappa(\rho)
\,
\left(2\left(\sqrt{\xi^*}+\sqrt{2}\right)\sigma_*
+
\varsigma^*_\zs{n}(r)
\right)
\frac{v_\zs{n}}{\sqrt{n}}
\,.
$$
Moreover, if the sequencies \eqref{Or.17} satisfy the properties
 \eqref{Or.25}--\eqref{Or.27}
then, for any $\gamma>0$,
$$
\lim_\zs{n\to\infty}
\sup_\zs{\varrho_\zs{n}\le \rho\le 1/6}
\cD_\zs{n}(\rho,r)/n^\gamma=0\,,
$$
where $\varrho_\zs{n}$ is any slowly decreasing sequence, i.e. satisfying
\eqref{Or.20}.
\end{theorem}
\begin{remark}\label{Re.3.4}
The inequality \eqref{Si.6} is used to prove the asymptotic efficiency
of the estimator \eqref{Or.15} (see, 
\cite{GaPe6}). 
To obtain the minimal
asymptotic quadratic risk, one has to take a weighted least squares estimator
\eqref{Or.8} with weights of type
\eqref{Or.23}, where the parameter $\alpha$ depends on unknown smoothness of unknown
function $S$. So one can't use this estimator directly. The oracle inequality allows us
to overcome this difficulty because the upper bound is majorized up to a multiplicative
and a additive constants by the minimal quadratic risk over all weighted estimators 
including the optimal one. Taking into account that the multiplicative constant tends to
unity, as $n\to\infty$, and the additive constant is negligible in comparison with any
degree of $n$ and, in particular, with the optimal convergence rate and then multiplying
the inequality \eqref{Si.6} by the optimal convergence rate, one obtains the asymptotically minimal upper bound for the quadratic risk of the estimator \eqref{Or.15}. The last result
means that this estimator is asymptotically efficient.
\end{remark}
\vspace*{5mm}
\section{Proofs}\label{Se.P}
\subsection{Proof of Theorem~\ref{Th.Or.1}}

First of all, note that we can represent
the empirical squared error  $\Er_\zs{n}(\lambda)$ by the following way
\begin{equation}\label{P.1}
\Er_\zs{n}(\lambda)=J_\zs{n}(\lambda)+
2\sum^n_\zs{j=1}\lambda(j)\theta'_\zs{j,n}+
\|S\|^2_\zs{n}-\rho\,
\wh{P}_\zs{n}(\lambda)
\end{equation}
with $\theta'_\zs{j,n}\,=\,\wt{\theta}_\zs{j,n}\,-\,\theta_\zs{j,n}\wh{\theta}_\zs{j,n}$.
The second term is most difficult to handle in the right-hand part of \eqref{P.1}. To estimate, we decompose it in the sum of terms and we apply appropriate technique to each of them. By setting
\begin{equation}\label{P.1-1}
\varsigma_\zs{j,n}=
\E_\zs{S}\,\xi^2_\zs{j,n}
=\frac{1}{n}\sum^n_\zs{l=1}\sigma^2_\zs{l}\phi^2_\zs{j}(x_l)
\quad
\mbox{and}
\quad
\mu_\zs{j,n}=\xi^2_\zs{j,n}-\varsigma_\zs{j,n}
\,,
\end{equation}
 we find that
\begin{equation}\label{P.1-2}
\theta'_\zs{j,n}=\frac{1}{\sqrt{n}}\theta_\zs{j,n}\xi_\zs{j,n}+
\frac{1}{n}\mu_\zs{j,n}
+\frac{1}{n}\,
(\varsigma_\zs{j,n}-\wh{\varsigma}_\zs{n})\,.
\end{equation}
Moreover, we can represent $\mu_\zs{j,n}$ as
\begin{equation}\label{P.2}
\mu_\zs{j,n}=\frac{1}{n}
\sum^n_\zs{l=1}\sigma^2_\zs{l}\phi^2_\zs{j}(x_l)
\eta_\zs{l}+
2\,\sum^n_\zs{l=2}\tau_\zs{j,l}\,\xi_\zs{l}
=\mu'_\zs{j,n}+2\mu''_\zs{j,n}
\end{equation}
with $\eta_\zs{l}=\xi^2_\zs{l}-1$ and
$$
\tau_\zs{j,l}=\frac{1}{n}\sigma_\zs{l}\phi_\zs{j}(x_l)
\sum^{l-1}_\zs{k=1}\sigma_\zs{k}\phi_\zs{j}(x_\zs{k})\,\xi_\zs{k}\,.
$$
Now we set
\begin{equation}\label{P.2-1}
N'(\lambda)=\sum^n_\zs{j=1}\lambda(j)\,
\mu'_\zs{j,n}
\quad\mbox{and}\quad
N''(\lambda)=\frac{1}{\sqrt{n\varsigma_\zs{n}}}\sum^n_\zs{j=1}\ov{\lambda}(j)\,
\mu''_\zs{j,n}\,\Chi_\zs{\{\varsigma_\zs{n}>0\}}\,,
\end{equation}
where $\ov{\lambda}(j)=\lambda(j)/|\lambda|$.
In the Appendix we show that
\begin{equation}\label{P.3}
\sup_\zs{\lambda\in\Lambda}\,\E_\zs{S}\, |N'(\lambda)|\,\le
\,
\sigma_*(
v_\zs{n}
+u_\zs{1,n})\frac{\sqrt{\xi^*}}{\sqrt{n}}
\end{equation}
and
\begin{equation}\label{P.4}
\sup_\zs{\lambda\in\bbr^n}\,\E_\zs{S}(N''(\lambda))^2\le\,\frac{2\sigma_*}{n}\,.
\end{equation}
Now, for any $\lambda\in\Lambda$, we rewrite \eqref{P.1} as
\begin{align*}
\Er_\zs{n}(\lambda)&=J_\zs{n}(\lambda)
+\frac{2}{n}N'(\lambda)
+4\sqrt{P_\zs{n}(\lambda)}
N''(\lambda)\\
&+2M(\lambda)
+\frac{2}{n}\,\Delta(\lambda)+\|S\|^2_\zs{n}
-\rho\wh{P}_\zs{n}(\lambda)\,,
\end{align*}
where $P_\zs{n}(\lambda)$ is defined in \eqref{Or.13},
\begin{align}\label{P.5}
\Delta(\lambda)=
\sum^n_\zs{j=1}\,\lambda(j)\,(\varsigma_\zs{j,n}-\wh{\varsigma}_\zs{n})
\quad\mbox{and}\quad
M(\lambda)=n^{-1/2}
\sum^n_\zs{j=1}\lambda(j)\theta_\zs{j,n}\xi_\zs{j,n}\,.
\end{align}
Further we estimate the term $\Delta(\lambda)$. Setting
\begin{equation}\label{P.5-1}
\ov{\varsigma}_\zs{j,n}=
\varsigma_\zs{j,n}-\varsigma_\zs{n}=
\frac{1}{n}\sum^n_\zs{l=1}\sigma^2_\zs{l}( \phi^2_\zs{j}(x_\zs{l})-1)\,,
\end{equation}
we obtain that
\begin{align}\nonumber
|\Delta(\lambda)|&\le
\,|\,\sum^n_\zs{j=1}\lambda(j) \ov{\varsigma}_\zs{j,n}\,|
+v_\zs{n}|\wh{\varsigma}_\zs{n}-\varsigma_\zs{n}|\\ \label{P.6}
&\le
\sigma_*u_\zs{1,n}+v_\zs{n}|\wh{\varsigma}_\zs{n}-\varsigma_\zs{n}|\,.
\end{align}
Now from \eqref{P.1}
we obtain that, for some
fixed $\lambda_0\in\Lambda$,
\begin{align*}
\Er_\zs{n}(\wh{\lambda})-\Er_\zs{n}(\lambda_0)&=J(\wh{\lambda})-J(\lambda_0)+
2M(\wh{\vartheta})
+\frac{2}{n}N'(\wh{\vartheta})\\
&+ 4\sqrt{P_\zs{n}(\wh{\lambda})}
N''(\wh{\lambda})
-4\sqrt{P_\zs{n}(\lambda_\zs{0})} N''(\lambda_0)
\\
&
-\rho \wh{P}_\zs{n}(\wh{\lambda})
+\rho \wh{P}_\zs{n}(\lambda_0)
+\frac{2}{n}
\left(
\Delta(\wh{\lambda})
-
\Delta(\lambda_\zs{0})
\right)\,,
\end{align*}
where $\wh{\vartheta}=\wh{\lambda}-\lambda_\zs{0}$.

By the  definition of $\wh{\lambda}$
in \eqref{Or.14} and by \eqref{P.6} we get
\begin{align*}
\Er_\zs{n}(\wh{\lambda})-\Er_\zs{n}(\lambda_0)&\le
2M(\wh{\vartheta})+
\frac{4\sigma_*u_\zs{1,n}+4v_\zs{n}|\wh{\varsigma}_\zs{n}-\varsigma_\zs{n}|}{n}\\
&+\frac{2}{n}N'(\wh{\vartheta})
+4 \sqrt{P_\zs{n}(\wh{\lambda})} N''(\wh{\lambda})-\rho \wh{P}_\zs{n}(\wh{\lambda})
\\
&
+\rho \wh{P}_\zs{n}(\lambda_0)
-4\sqrt{P_\zs{n}(\lambda_0)}
N''(\lambda_0)\,.
\end{align*}
Moreover, making use of  the inequality
\begin{equation}\label{P.6-1}
2|ab|\le \epsilon a^2+\epsilon^{-1}b^2
\end{equation}
with $\epsilon=\rho/4$
and taking into account
 the definition of the penalty term in \eqref{Or.12} we deduce, for any $\lambda\in\Lambda$,
\begin{align*}
4 \sqrt{P_\zs{n}(\lambda)}
|N''(\lambda)|&\le \rho P_\zs{n}(\lambda)
+4\frac{(N''(\lambda))^2}{\rho}\\
&\le \rho \wh{P}_\zs{n}(\lambda)
+\rho \frac{|\lambda|^2|\wh{\varsigma}_\zs{n}-\varsigma_\zs{n}|}{n}+\frac{4(N''(\lambda))^2}{\rho}\,.
\end{align*}
Thus from here it follows that
\begin{align}\label{P.8}
\Er_\zs{n}(\wh{\lambda})&\le\Er_\zs{n}(\lambda_0)+
2M(\wh{\vartheta})+\Upsilon_\zs{n}
+2\rho \wh{P}_\zs{n}(\lambda_0)\,,
\end{align}
where
$$
\Upsilon_\zs{n}
=\frac{4}{n}N^*_\zs{1} +\frac{8}{\rho }(N^*_\zs{2})^2
+\frac{4\sigma_*u_\zs{1,n}}{n}
+\frac{4+2\rho}{n}v_\zs{n}|\wh{\varsigma}_\zs{n}-\varsigma_\zs{n}|
$$
with $N^*_\zs{1}=\sup_\zs{\lambda\in\Lambda}|N'(\lambda)|$ and
$N^*_\zs{2}=\sup_\zs{\lambda\in\Lambda}|N''(\lambda)|$.
Moreover, note that the bounds \eqref{P.3}, \eqref{P.4} and \eqref{P.6} imply
that
\begin{align}\label{P.9}
\E_\zs{S}\Upsilon_\zs{n}
\le \Upsilon^*_\zs{n}(\rho)
\frac{\sigma_*}{n}+\frac{4+2\rho}{n}v_\zs{n}
\E_\zs{S}|\wh{\varsigma}_\zs{n}-\varsigma_\zs{n}|\,,
\end{align}
where  the function $\Upsilon^*_\zs{n}(\rho)$ is defined in \eqref{Or.16}.

Now we study the second term in \eqref{P.5}. First, note that
for any nonrandom vector $\vartheta=(\vartheta(1),\ldots,\vartheta(n))'\in\bbr^n$
 Lemma~\ref{Le.A.5}
implies 
\begin{equation}\label{P.12}
\E_\zs{S} M^2(\vartheta)
\le
\frac{\sigma_*}{n}\,
\sum^n_\zs{j=1}\vartheta^2(j)
\theta^2_\zs{j,n}
=\sigma_*
\frac{\|S_\zs{\vartheta}\|^2_n}{n}\,,
\end{equation}
where 
$$
S_\zs{\vartheta}=\sum^n_\zs{j=1}\vartheta(j)\theta_\zs{j,n}\phi_\zs{j}\,.
$$
We set now
$$
Z^*=\sup_\zs{\vartheta\in \Lambda_\zs{1}}\frac{n M^2(\vartheta)}{\|S_\zs{\vartheta}\|^2_n}
\quad\mbox{with}\quad
\Lambda_\zs{1}=\Lambda-\lambda_\zs{0}\,.
$$
To estimate this term we apply the inequality \eqref{P.12}, i.e. 
\begin{equation}\label{P.14}
\E_\zs{S}\,Z^*\,\le\,
\sum_\zs{\vartheta\in \Lambda_\zs{1}}\frac{n\E_\zs{S}\,M^2(\vartheta)}{\|S_\zs{\vartheta}\|^2_n}
\le \nu\sigma_*\,.
\end{equation}
Moreover, making use of inequality \eqref{P.6-1} with  $\epsilon=\rho\|S_\zs{\vartheta}\|_\zs{n}$,
we get
\begin{equation}\label{P.13}
2|M(\vartheta)|\le\rho \|S_\zs{\vartheta}\|^2_\zs{n}+
\frac{Z^*}{n\rho}\,.
\end{equation}
Now we estimate  $\|S_\zs{\vartheta}\|^2_\zs{n}$. We have
\begin{equation}\label{P.15}
\|S_\zs{\vartheta}\|^2_\zs{n}-
\|\wh{S}_\zs{\vartheta}\|^2_\zs{n}=
\sum^n_\zs{j=1}\,\vartheta^2(j)
(\theta^2_\zs{j,n}-\wh{\theta}^2_\zs{j,n})\le -2 M_\zs{1}(\vartheta)
\end{equation}
with
$$
M_\zs{1}(\vartheta)=
\frac{1}{\sqrt{n}}\sum^n_\zs{j=1}\vartheta^2(j)\theta_\zs{j,n}\xi_\zs{j,n}\,.
$$
Now, taking into account that $|\vartheta(j)|\le 1$
 for any $\vartheta\in\Lambda_\zs{1}$, we obtain
$$
\E_\zs{S} M^2_\zs{1}(\vartheta)\le \sigma_*
\frac{\|S_\zs{\vartheta}\|^2_n}{n}\,.
$$
Putting
$$
Z^*_\zs{1}=\sup_\zs{\vartheta\in \Lambda_\zs{1}}
\frac{nM^2_\zs{1}(\vartheta)}{\|S_\zs{\vartheta}\|^2_n}\,,
$$
 we get
\begin{equation}\label{P.16}
\E_\zs{S}\,Z^*_\zs{1}\,\le \nu\sigma_*\,.
\end{equation}
Therefore, applying inequality \eqref{P.13} for $M_\zs{1}(\vartheta)$ in 
\eqref{P.15} we deduce the upper bound for $\|S_\zs{\vartheta}\|^2_\zs{n}$, i.e.
\begin{align}\label{P.17}
\|S_\zs{\vartheta}\|^2_\zs{n}&\le\,\frac{1}{1-\rho}
\|\wh{S}_\zs{\vartheta}\|^2_\zs{n}
+\,
\frac{Z^*_\zs{1}}{n\rho(1-\rho)}\,.
\end{align}
Taking into account this inequality in \eqref{P.13} we obtain that
\begin{align*}
2M(\vartheta)&\le\,\frac{\rho }{1-\rho }
\|\wh{S}_\zs{\vartheta}\|^2_\zs{n}+
\frac{Z^*+Z^*_\zs{1}}{n\rho(1-\rho)}\\
&\le
\frac{2\rho(\Er_n(\lambda)+\Er_n(\lambda_0))}{1-\rho}
+
\frac{Z^*+Z^*_\zs{1}}{n\rho(1-\rho)}\,.
\end{align*}
Therefore  \eqref{P.8} implies that
\begin{align*}
\Er_\zs{n}(\wh{\lambda})&
\le \frac{1+\rho }{1-3\rho }
\Er_\zs{n}(\lambda_0)+
\frac{1-\rho }{1-3\rho}\Upsilon_\zs{n}\\
&+
\frac{Z^*+Z^*_\zs{1}}{n\rho(1-3\rho )}
+\frac{2\rho(1-\rho)}{1-3\rho}\wh{P}_\zs{n}(\lambda_\zs{0})\,,
\end{align*}
Now by  inequalities
 \eqref{P.14}--\eqref{P.16}
we get that
\begin{align*}
\E_\zs{S}\Er_\zs{n}(\wh{\lambda})&
\le \frac{1+\rho}{1-3\rho}
\E_\zs{S}\Er_\zs{n}(\lambda_\zs{0})+
\frac{1-\rho}{1-3\rho}\E_\zs{S}\,\Upsilon_\zs{n}\\
&
+\frac{2\nu\sigma_*}{n\rho(1-3\rho)}
+\frac{2\rho(1-\rho)}{1-3\rho}\,
\E_\zs{S} \wh{P}_\zs{n}(\lambda_\zs{0})\,.
\end{align*}
By making use of inequality \eqref{P.9} and Lemma~\ref{Le.A.1}
we come to Theorem~\ref{Th.Or.1}.
\endproof

\subsection{Proof of Lemma~\ref{Le.Si.1}}
First notice that from \eqref{Or.7} we obtain that
\begin{align*}
\wh{\varsigma}_\zs{n}-\varsigma_\zs{n}&\,=\,
\sum^n_\zs{j=d_\zs{n}+1}\theta^2_\zs{j,n}+
\frac{2}{\sqrt{n}}\sum^n_\zs{j=d_\zs{n}+1}\,\,\,
\theta_\zs{j,n}\,\xi_\zs{j,n}\\ 
&+\,
n^{-1}\,\sum^n_\zs{j=d_\zs{n}+1}\,
\mu_\zs{j,n}\,+\,
n^{-1}\,\sum^n_\zs{j=d_\zs{n}+1}\,\,\ov{\varsigma}_\zs{j,n}\,-\,
\frac{d_\zs{n}}{n}\,\varsigma_\zs{n}\\
&:=\Delta_\zs{1}+\frac{2}{\sqrt{n}}\Delta_\zs{2}+\frac{1}{n}\Delta_\zs{3}+
\frac{1}{n}
\Delta_\zs{4}-\frac{d_\zs{n}}{n}\,\varsigma_\zs{n}\,,
\end{align*}
where $\mu_\zs{j,n}$ and $\ov{\varsigma}_\zs{j,n}$ are defined in
\eqref{P.1-2} and \eqref{P.5-1} respectively.

We estimate the first term by Lemma~\ref{Le.A.4} for $S\in W^1_\zs{r}$. We have
$$
\Delta_\zs{1}\le \frac{4r}{d^2_\zs{n}}\,.
$$
The next term we estimate with the help of Lemma~\ref{Le.A.5}.
We get that
$$
\E_\zs{S}(\Delta_\zs{2})^2\le \sigma_*\Delta_\zs{1}
\le \sigma_*\frac{4r}{d^2_\zs{n}}\,.
$$
By \eqref{P.2} and \eqref{P.2-1} we can represent $\Delta_\zs{3}$ as
$$
\Delta_\zs{3}=N'(\lambda_\zs{I})
+2|\lambda_\zs{I}|\sqrt{n\varsigma_\zs{n}}N''(\lambda_\zs{I})
$$
with the vector 
$\lambda_\zs{I}=(\lambda_\zs{I}(1)\,,\ldots,\,\lambda_\zs{I}(n))'$ 
having the indicator components, i.e. 
$\lambda_\zs{I}(j)=\Chi_\zs{\{j>d_\zs{n}\}}$. By estimating in \eqref{A.0} 
$\phi^2_\zs{j}$ by $2$ we obtain
$$
\E_\zs{S}\,|N'(\lambda_\zs{I})|\le\, 2\sigma_*\,\sqrt{\xi^*}\,\sqrt{n}\,.
$$ 
Thus the upper bound
\eqref{P.4} implies
\begin{align*}
\E_\zs{S}|\Delta_\zs{3}|\le 2\sigma_*(\sqrt{\xi^*}+\sqrt{2})\sqrt{n}=\ov{\sigma}\sqrt{n}\,.
\end{align*}

Moreover, due to Lemma~\ref{Le.A.3} with $k=0$, one has
\begin{align*}
|\Delta_\zs{4}|\,&=\,
\left|n^{-1}\,\sum^n_\zs{l=1}\,\sigma^2_\zs{l}\,\sum^{n}_\zs{j=d_\zs{n}+1}
( \phi^2_\zs{j}(x_\zs{l})-1)
\right|\\
&\le
\frac{\sigma_*}{n}
\sum^n_\zs{l=1}\,\left|\,\sum^{n}_\zs{j=1}\,( \phi^2_\zs{j}(x_\zs{l})-1)\right|+
\frac{\sigma_*}{n}
\sum^n_\zs{l=1}\,\left|\,\sum^{d_\zs{n}}_\zs{j=1}\,( \phi^2_\zs{j}(x_\zs{l})-1)\right|
\le 2\sigma_*\,.
\end{align*}
Hence 
Lemma~\ref{Le.Si.1}.
\endproof
\medskip
\medskip

\section{Numerical example}\label{Se.Co}

In this Section, we compare via simulations 
the adaptive procedure proposed for the heteroscedastic regression in
\cite{EfPi} (section 4.1) with that of \eqref{Or.15}.

 Consider the  model \eqref{I.1} with
$$
S(x)\,=\,x\, \sin(2\pi x)\,+\,x^2\,(1-x)\,\cos(4\pi x)
\quad\mbox{and}\quad
\sigma^2_\zs{j}=1+S^2(x_\zs{j})\,,
$$
 assuming that
$(\xi_\zs{k})_\zs{k\ge 1}$ follow the gaussian distribution with zero
mean and unit variance. 

In the procedure \eqref{Or.15} we take the weight vectors defined in
\eqref{Or.23} with $k^*=100+\sqrt{\ln n}$, $\ve=1/\ln n$, $m=\ln^2 n$,
$\rho=1/(3+\ln^2n)$
and 
$$
\omega_\zs{\beta}=10+(A_\zs{\beta}tn)^{1/(2\beta+1)}\,.
$$
Moreover, we make use of the estimate \eqref{Si.4} with $d_\zs{n}=n^{1/3}$.
\medskip
\medskip

 The results of
simulations are reported in Table 1.
\begin{center}
                {\it Table 1}\\
\vspace*{7mm}
\begin{tabular}{|c|c|c|c|c|c|}
\hline $\hat{E}\,\|S^*-S\|^2_\zs{n}$
&$\hat{E}\,\|\tilde{S}-S\|^2_\zs{n}$ &$n$\\[2mm]
\hline 0.260&0.410& 21\\
\hline 0.148&0.427& 41\\
\hline 0.058&0.476& 101\\
\hline 0.034&0.430& 201\\
\hline 0.019&0.448& 401\\
\hline
\end{tabular}
\end{center}

\vspace*{10mm}

The columns of Table 1 with the headings 
$\hat{E}\,\|S^*-S\|^2_\zs{n}$, $\hat{E}\,\|\tilde{S}-S\|^2_\zs{n}$
and $n$
report, respectively, the empirical quadratic risk
for the procedure \eqref{Or.23},  
the empirical quadratic risk
for the procedure from \cite{EfPi}
and
the sample size. To calculate the empirical risks  50 independent Monte Carlo simulations
were performed.

Table 1 shows that in comparison with the procedure from \cite{EfPi} our adaptive estimator performs resonably well for the small sample sizes.

\medskip

\renewcommand{\theequation}{A.\arabic{equation}}
\renewcommand{\thetheorem}{A.\arabic{theorem}}
\renewcommand{\thesubsection}{A.\arabic{subsection}}
\section{Appendix}\label{Se.A}
\setcounter{equation}{0}
\setcounter{theorem}{0}

\subsection{Proof of \eqref{P.3}}\label{Su.A.1}
First note that we can represent the term $N'(\lambda)$
as
$$
N'(\lambda)=\sum^n_\zs{l=1}\,g_\zs{l,n}\,\eta_\zs{l}
\quad
\mbox{with}
\quad
g_\zs{l,n}=\frac{\sigma^2_\zs{l}}{n}\sum^n_\zs{j=1}\lambda(j)\phi^2_\zs{j}(x_\zs{l})\,.
$$
Recalling that $\E\,\eta^2_\zs{1}=\xi^*-1$ we calculate
$$
\E_\zs{S}\,(N'(\lambda))^2=
\frac{\xi^*-1}{n^2}\sum^n_\zs{l=1}\E_\zs{S}\sigma^4_\zs{l}
\left(\sum^n_\zs{j=1}\lambda(j)\phi^2_\zs{j}(x_\zs{l})\right)^2\,.
$$
Therefore for any vector $\lambda\in\bbr^n$
\begin{equation}\label{A.0}
\E_\zs{S}\,|N'(\lambda)|\le \sigma_*
\frac{\sqrt{\xi^*}}{\sqrt{n}}\,\max_\zs{1\le l\le n}
\,|\sum^n_\zs{j=1}\lambda(j)\phi^2_\zs{j}(x_\zs{l})\,|\,.
\end{equation}
Thus taking into account here definitions \eqref{Or.17}
we come to inequality \eqref{P.3}.
\endproof

\subsection{Proof of \eqref{P.4}}\label{Su.A.2}
By putting $g_\zs{l}=\sum^n_\zs{j=1}\ov{\lambda}(j)\tau_\zs{j,l}$
and taking into account that the random variables $(\xi_\zs{k})_\zs{1\le k\le n}$
are independent of $(\sigma_\zs{k})_\zs{1\le k\le n}$
we obtain that
\begin{equation}\label{A.8}
\E_\zs{S}\left((N''(\lambda))^2\,|\,\sigma_\zs{1},\ldots,\sigma_\zs{n}\right)
=
\Chi_\zs{\{\varsigma_n>0\}}
\left(\sum^n_\zs{k=1}\sigma^2_\zs{k}\right)^{-1}\sum^n_\zs{l=2}\wh{g}_\zs{l}\,,
\end{equation}
where
\begin{align*}
\wh{g}_\zs{l}=\E(g^2_\zs{l}\,|\,\sigma_\zs{1},\ldots,\sigma_\zs{n})
=
\frac{\sigma^2_\zs{l}}{n^2}\sum^{l-1}_\zs{k=1}\sigma^2_\zs{k}
\left(\sum^n_\zs{j=1}\ov{\lambda}(j)\phi_\zs{j}(x_\zs{l}) \phi_\zs{j}(x_\zs{k})\right)^2\,.
\end{align*}
Therefore the orthonormality property \eqref{Or.5} 
implies that for any $\lambda\in\bbr^n$
\begin{align*}
\wh{g}_\zs{l}&\le \sigma_*
\frac{\sigma^2_\zs{l}}{n^2}\sum^{n}_\zs{k=1}
\left(\sum^n_\zs{j=1}\ov{\lambda}(j)\phi_\zs{j}(x_\zs{l}) \phi_\zs{j}(x_\zs{k})\right)^2\\
&=
\sigma_*\frac{\sigma^2_\zs{l}}{n}\,
\,\sum^n_\zs{j=1}\ov{\lambda}^2(j)\phi^2_\zs{j}(x_\zs{l})\,\le\,
\frac{2\sigma_*}{n}\,\sigma^2_\zs{l}\,.
\end{align*}
Now by making use of this inequality in \eqref{A.8} we get
\eqref{P.4}.

\endproof
\medskip

\subsection{Technical lemma}\label{Su.A.3}
\begin{lemma}\label{Le.A.1}
For any $n\ge 1$ and $\lambda\in\Lambda$,
\begin{align*}
\E_\zs{S}\wh{P}_\zs{n}(\lambda)
\le \E_\zs{S}\,\Er_\zs{n}(\lambda)+\frac{v_\zs{n}}{n}
\E_\zs{S}|\wh{\varsigma}_\zs{n}-\varsigma_\zs{n}|+
\frac{\sigma_*u_\zs{2,n}}{n}\,.
\end{align*}
\end{lemma}
\proof
Indeed, by the definition of $\Er_\zs{n}(\lambda)$
 we have
\begin{align*}
\Er_\zs{n}(\lambda)=
\sum^n_\zs{j=1}
\left((\lambda(j)-1)\theta_\zs{j,n}+\lambda(j)
\frac{1}{\sqrt{n}}\xi_\zs{j,n}\right)^2\,.
\end{align*}
Therefore,
\begin{align*}
\E_\zs{S} \Er_\zs{n}(\lambda)\ge
\E_\zs{S}\frac{1}{n}\sum^n_\zs{j=1}\lambda^2(j)\,\xi^2_\zs{j,n}
=\E_\zs{S}\frac{1}{n}\sum^n_\zs{j=1}\lambda^2(j)\,\varsigma_\zs{j,n}\,,
\end{align*}
where the sequence $(\varsigma_\zs{j,n})$ is defined in \eqref{P.1-1}. 
Moreover, note that the last term can be estimated as
$$
\left|\sum^n_\zs{j=1}\lambda^2(j)\varsigma_\zs{j,n}-|\lambda|^2\varsigma_\zs{n}\right|=
\left|\frac{1}{n}\sum^n_\zs{l=1}\sigma^2_\zs{l}\,
\sum^n_\zs{j=1}\lambda^2(j)\,(\phi^2_\zs{j}(x_\zs{l})-1)\right|
\le \sigma_*u_\zs{2,n}\,.
$$
We recall that the definition of the set $\Lambda$ and
the definition of $v_\zs{n}$ in \eqref{Or.17}
imply that $|\lambda|^2\le v_\zs{n}$
 for $\lambda\in\Lambda$. Therefore for any $\lambda\in\Lambda$
\begin{align*}
\sum^n_\zs{j=1}\lambda^2(j)\,\varsigma_\zs{j,n}
&\ge
|\lambda|^2 \wh{\varsigma}_\zs{n}-\sigma_*u_\zs{2,n}-|\lambda|^2
 |\wh{\varsigma}_\zs{n}-\varsigma_\zs{n}|\\
&\ge
|\lambda|^2 \wh{\varsigma}_\zs{n}-\sigma_*u_\zs{2,n}-v_\zs{n}
 |\wh{\varsigma}_\zs{n}-\varsigma_\zs{n}|\,.
\end{align*}
Hence
the desired inequality.
\endproof

\subsection{ Properties of trigonometric basis}\label{Su.A.4}

\begin{lemma}\label{Le.A.3}
For any $k\ge 0$,
\begin{equation}\label{A.2}
\sup_{N\ge 2}\quad\sup_{x\in [0,1]}N^{-k}\left|\sum_{l=2}^N\,l^k\,
(\phi^2_\zs{l}(x)-1)\,
\right|\le 2^k\,.
\end{equation}
\end{lemma}
\noindent{\bf Proof.}
Due to the properties of the trigonometric functions, we get
\begin{align*}
\sum_{l=2}^N l^k
\,(\phi^2_\zs{l}(x)-1)\,
&=\sum_{1\le l\le N/2}(2l)^k\,\cos(4\pi lx)\\
&-\sum_{1\le l\le (N-1)/2}(2l+1)^k\,\cos(4\pi lx)\,.
\end{align*}
This yields
\begin{align*}
\left|\sum_{l=2}^N\,l^k\,(\phi^2_l(x)-1)\right|&
\le\left|\sum_{1\le l\le (N-1)/2}\left((2l+1)^k-(2l)^k\right)\cos(4\pi lx)\right|+N^{k}\\
&\le \sum_{1\le l\le (N-1)/2} \left((2l+1)^k-(2l)^k\right)+N^{k}\\
&=\sum_{1\le l\le (N-1)/2}\sum_{j=0}^{k-1}\left(\stackrel{k}{j}\right)(2l)^j
+N^{k}\,.
\end{align*}
This implies (\ref{A.2}).

\endproof

\begin{lemma}\label{Le.A.4} For any function $S\in  W_r^k$,
\begin{equation}\label{A.4}
\sup_{n\ge 1}\sup_{1\le m\le n-1}\,m^{2k}\,
\left(\sum_{j=m+1}^{n}\,\theta_\zs{j,n}^2
\right)
\,\le\,\frac{4r}{\pi^{2(k-1)}}\,.
\end{equation}
\end{lemma}
\noindent{\bf Proof.} First, note that any function
$S$ from  $W_r^k$ can be represented by its Fourier series, i.e.
$S=\sum^\infty_{j=1}\theta_j\phi_j$ with the coefficients defined
by \eqref{Si.3-1}.
By denoting the residual term
for $S$ as
$$
\Delta_m(x)\,=\,S\,-\,\sum^m_{j=1}\theta_j\phi_j
\,=\,\sum_{j=m+1}^{\infty}\theta_\zs{j}\phi_j(x)\,,
$$
 we obtain
that
\begin{align*}
\sum_{j=m+1}^{n}\,\theta_\zs{j,n}^2
\,=\,
\inf_{\alpha_1,\ldots,\alpha_m}\|S-\sum_{j=1}^{m}\,\alpha_j\,\phi_j\|^2_n
\,\le\,\|\Delta_m\|^2_n\,.
\end{align*}
Moreover, it is easy to deduce that
\begin{align*}
\|\Delta_m\|^2_n&\,=\,n^{-1}\sum_{k=1}^n\,\Delta^2_m(x_k)
=\sum_{k=1}^n\int_{x_\zs{k-1}}^{x_\zs{k}}\,\Delta^2_m(x_k)\d x\\
&\le\, 2\int_0^1\,\Delta^2_m(x)\d x\,+\,2\sum_{k=1}^n\int_{x_\zs{k-1}}^{x_\zs{k}}\,
(\Delta_m(x_k)-\Delta_m(x))^2\d x\,.
\end{align*}
The last term in this inequality we estimate as
\begin{align*}
(\Delta_m(x_k)-\Delta_m(x))^2&=
\left(\int_{x}^{x_\zs{k}}\,\dot{\Delta}_m(z)\d z\right)^2\\
&\le n^{-1}\int_{x_\zs{k-1}}^{x_\zs{k}}\,(\dot{\Delta}_m(z))^2\d z\,.
\end{align*}
Therefore,
\begin{align*}
\|\Delta_m\|_n^2&\le 2\|\Delta_m\|^2+\frac{2}{n^2}\|\dot{\Delta}_m\|^2\\
&=2\sum_{j=m+1}^{\infty}\theta_\zs{j}^2+
\frac{2}{n^2}\sum_{j=m+1}^{\infty}\theta_\zs{j}^2
\|\dot{\phi}_\zs{j}\|^2.
\end{align*}
Now note that by the representation of the set $W^k_\zs{r}$
 in the form \eqref{Si.3} we can estimate the first term in the last 
inequality as
\begin{align*}
\sum_{j=m+1}^{\infty}\theta_\zs{j}^2=
\sum_{j=m+1}^{\infty}\theta_\zs{j}^2\frac{a_\zs{j}}{a_\zs{j}}
\le \frac{r}{a_\zs{m+1}}\le  \frac{r}{(\pi m)^{2k}}\,.
\end{align*}
Similarly, we find that
\begin{align*}
\sum_{j=m+1}^{\infty}\theta_\zs{j}^2\|\dot{\phi}_\zs{j}\|^2
\le \sup_\zs{j\ge m+1}\frac{\|\dot{\phi}_\zs{j}\|^2}{a_\zs{j}}r
\le \sup_\zs{j\ge m+1}\frac{\|\dot{\phi}_\zs{j}\|^2}{\|\phi^{(k)}_\zs{j}\|^2}r
\le \frac{r}{(\pi m)^{2(k-1)}}\,.
\end{align*}
Therefore, for $m\le n$ we get that
\begin{align*}
\frac{1}{n^2}\sum_{j=m+1}^{\infty}\theta_\zs{j}^2\|\dot{\phi}_\zs{j}\|^2
\le \frac{r}{\pi^{2(k-1)} m^{2k}}\,.
\end{align*}
This implies \eqref{A.4}.
\endproof

\begin{lemma}\label{Le.A.5}
Let $\xi_\zs{j,n}$ be defined in \eqref{Or.7} for the model \eqref{I.1}.
Then, for any real numbers $f_\zs{1},\ldots,f_\zs{n}$,
\begin{equation}\label{A.5}
\E\,\left(\sum_{j=1}^{n}\,f_\zs{j}\,\xi_\zs{j,n}\right)^2\le\,
\sigma_*\sum_{j=1}^{n}\,f^2_\zs{j}\,.
\end{equation}
\end{lemma}
\noindent{\bf Proof.}
Due to the definition of  $\xi_\zs{j,n}$, one has
$$
\sum_{j=1}^{n}\,f^2_\zs{j}\,\xi_\zs{j,n}=
\sum_{l=1}^{n}\,\sigma_l\,\wt{f}_\zs{l}\,\xi_\zs{l}
\quad
\mbox{with}
\quad
\wt{f}_\zs{l}=\frac{1}{\sqrt{n}}\,\sum^n_\zs{j=1}\,f_\zs{j}\,\phi_j(x_l)\,.
$$
 Moreover
\begin{align*}
\E\left(\sum_{j=1}^{n}\,f_\zs{j}\,\xi_\zs{j,n}\right)^2&=
\,\sum_{l=1}^{n}\,\sigma^2_\zs{l}\wt{f}^2_\zs{l}\,
\le\,
\sigma_*\,\sum_{l=1}^{n}\,\wt{f}^2_\zs{l}\\
&=
\sigma_*\sum_{i,j=1}^{n}\, f_\zs{i}\, f_\zs{j}\, (\phi_i,\phi_j)_n\,.
\end{align*}
The orthogonality  of the basis $(\phi_\zs{j})$
implies inequality  \eqref{A.5}.
Hence Lemma~\ref{Le.A.5}.
\endproof

\medskip

{\bf Acknowledgements}

The authors are grateful to the anonymous Reviewers for a careful reading of manuscript and
helpful comments which led to an improved presentation of the paper.

\end{document}